\documentclass[12pt]{article}%
\usepackage{amsfonts}
\usepackage{amssymb}
\usepackage{graphicx}
\usepackage{eurosym}

\usepackage{amsmath}%
\newcommand{\bx}{{\bf x}}
\newcommand{\by}{{\bf y}}

\begin{document}

\title{\vspace{-1cm} Optimisation via Slice Sampling}

\author{John R. Birge\\Booth School of Business \and Nicholas G. Polson\\Booth School of Business\footnote{{\em Address for correspondence:}
        John Birge, 5807 South Woodlawn Avenue, Chicago, IL 60637, U.S.A.
        E-mails: john.birge@chicagobooth.edu; ngp@chicagobooth.edu. The authors' work was supported by the University of Chicago Booth School of Business.}
        }

\date{First Draft: January 2012\\
This Draft: June 2012}

\maketitle
\begin{abstract}
\noindent In this paper, we develop a simulation-based approach to optimisation with multi-modal functions using slice sampling.
Our method specifies the objective function as an energy potential in a Boltzmann distribution and then we use
auxiliary exponential slice variables to provide samples for a variety of energy levels.
Our slice sampler draws uniformly over the augmented slice region. We identify the global modes by projecting the path of the chain back to the underlying space.
Four standard test functions are used to illustrate the methodology: Rosenbrock,
Himmelblau, Rastrigin,  and Shubert.
These functions demonstrate the flexibility of our approach as they include functions with  long
ridges (Rosenbrock), multi-modality (Himmelblau, Shubert) and many local modes dominated by one global (Rastrigin).
The methods described here are implemented in the {\tt R} package {\tt McmcOpt}.
\end{abstract}

\vspace{1in} {\bf Keywords:} Himmelblau, Rastrigin, Rosenbrock, Shubert, Boltzmann Distribution, Slice Sampling, Simulation, Stochastic Optimisation, MCMC, Markov chain.

\newpage

\section{Introduction}

Multi-modal objective functions pose difficulties for local search and derivative based methods.
Our simulation-based approach exploits a well-known duality between functional optimisation and sampling to find the mode of the Boltzmann distribution
with an energy potential specified by the objective function of interest.
We exploit auxiliary slice variables to augment the state space and develop a Markov chain which samples uniformly on the slice region whilst
traversing the desired modes in the original space. A simulation-based approach has the advantage of being derivative free and thus avoids some
of the problems associated with optimisation of multi-modal functions.
To illustrate our methodology, we test four standard global functions from the optimisation literature:
the Rosenbrock, Himmelblau, Rastrigin, and Shubert.

Our approach builds on the seminal work of
Pincus (1968, 1970), Geman and Geman (1985), and Geman (1990) who proposed simulation-based Metropolis-Hastings (MH) algorithms for Boltzmann distributions
in the context of constrained functional optimisation.
We use a slice sampling approach for additive functionals (Edwards and Sokal, 1988) that scales to high dimensions.
The Boltzmann distribution contains an energy level parameter which we can use to perform a sensitivity analysis.
Subsequent research has shown that  direct MCMC methods can be fraught with convergence difficulties as the 
associated chain can easily become stuck in a local mode;
consequentially, careful tuning of algorithms is generally required.

One of the advantages of our approach
is that it does not require additional tuning. For the algorithm to be efficient the chain has to have good mixing properties; 
see, Tweedie and Mengersen (1994), Polson (1996) and
Roberts and Rosenthal (1999), Mira and Tierney (2002) for theoretical convergence results on slice sampling.
From a practical perspective, in the four examples considered here, our slice sampler has remarkably good mixing properties as it samples uniformly over the higher dimensional
slice region. In all cases, with a reasonable set of energy levels, we can traverse the objective functions of interest within thousands of MCMC draws.

Other popular simulation-based methods range
from simulated annealing (Kirkpatrick et al, 1983, Geman, 1990),
direct and evolutionary Metropolis MCMC (Liu, et al, 2000), particle swarm (Kennedy and Eberhart, 1995), multi-set sampling (Leman, Chen and Lavine, 1999),
and stochastic guided search (Gramacy and Taddy, 2010, Gramacy and Lee, 2010, Taddy, Lee, Gray and Griffin, 2008). For example,
Janson and Middendorf (2005) illustrate particle swarm methods on the
Rosenbrock and Rastigin functions. Many of these methods require careful tuning. Our approach relies solely on the energy level as a free tuning parameter and
requires only a sensitivity analysis to this parameter. We follow the simulated tempering literature by focusing on a set
of pre-determined energy values although our method extends easily to the case of stochastic levels as used in the Wang Landau algorithm 

The rest of the paper is as follows. Section 2 describes our simulation-based optimisation procedure. Section 3 develops slice samplers for
the  Rosenbrock, Himmelblau, Rastrigin, and Shubert functions.
We show how to calculate the slice-set for each of the functions in turn and develop the associated MCMC algorithms.
Finally, Section 4 concludes with directions for future research.

\section{Simulation-based Optimisation via Slice Sampling}

The general problem we address is to find the set of minima
$ {\rm arg min}_{ \bx \in \mathcal{X}} f(\bx)$, for some domain $ \mathcal{X} $ for a given objective $f(\bx)$. 
This is the standard optimisation problem. Our method distinguishes itself from others in that we allow for
the function $f(\bx)$ to be multi-modal. We define the set of minima by
$$
\mathcal{X}_{min} = \{  \bx \in \mathcal{X} : f( \bx) = \min_\by f(\by) \} \; .
$$
We will find $\mathcal{X}_{min}  $ by simulation, however, we do not directly simulate $f(\bx)$, rather we exploit a well-known duality between optimisation and
finding the modes of the  Boltzmann
distribution with energy potential $f(\bx)$ defined by the density
$$
\pi_\kappa (\bx) = \exp \left \{ - \kappa f(\bx) \right \} / Z_\kappa \; \; {\rm for} \; \; \bx \in \mathcal{X} 
$$
where $ Z_\kappa = \int_{\mathcal{X}}  \exp \left \{ - \kappa f(\bx) \right \} d \bx $ is an appropriate
normalisation constant or partition function. Clearly the minima of $f(\bx)$ correspond to the modes of $\pi_\kappa(\bx)$.
One advantage of our simulation approach is that we do not require explicit knowledge of $Z_\kappa $.
The only tuning parameter is $\kappa$ which is an energy level parameter of the Boltzmann distribution.

The limiting cases of the Boltzmann distribution where  $ \kappa \rightarrow 0 $ or $ \kappa \rightarrow  \infty  $ 
are of particular interest. They both lead to a uniform measure but on different sets.
When $ \kappa \rightarrow 0 $, the limiting distribution, denoted by $ \pi_0 (\bx) $, is a uniform measure on the set $ \mathcal{X} $. 
When $ \kappa \rightarrow \infty $, the limiting distribution, denoted by $ \pi_\infty (\bx) $, is a uniform measure on the set of modes,
$ \mathcal{X}_{min} $. Therefore, if we can sample from the Boltzmann distribution we can identify the minima of the original function.
Specifically, we have
$$
\lim_{ \kappa \rightarrow \infty } \pi_\kappa (\bx) = \pi_{\infty} ( \bx ) = | \mathcal{X}_{min} |^{-1} \delta_{ \mathcal{X}_{min} } (\bx)
$$
where $ \delta $ denotes a Dirac measure.

Once our problem is re-written in terms of the Boltzmann distribution we can extend existing methods for finding the modes. For example,
Pincus (1968, 1970) proposed a Metropolis algorithm to simulate draws. 
We denote the realisation of the  Markov chain by $ X^{(0)} , X^{(1)} , \ldots , X^{(G)} , \ldots $
which has equilibrium distribution, $ \pi_\kappa (\bx)$. Then, under mild Harris recurrence conditions,
given any starting point $X^{(0)}$ and energy level $\kappa$, we have the limit
$$ 
\lim_{ G \rightarrow \infty} \mathbb{P} \left ( X^{(G)} \in A | X^{(0)}= \bx \right ) = \pi_\kappa (A )
$$
for any Borel sets $A$.
See Tierney (1994) and Azencott (1988) for further discussion. 
We can then use the ergodic mean $ \frac{1}{G} \sum_{g=1}^G X^{(g)} $ along the chain as an estimate of the mean and hence, 
in the uni-modal case as $ \kappa \rightarrow \infty $ this will find the mode.

There are, however, many possible Markov transition dynamics that have the appropriate equilibrium distribution.
Thus the main issue becomes which Markov chain to use. We argue for the use of slice sampling methods.
The practical insight of using augmentation and slice sampling is that
essentially we have put some volume back into the spiky multi-mode regions. After the chain has converged in the higher dimensional set,  we can
project the draws back down into the dimension of interest and the chain will have no difficulty in traversing the modes
even for lower energy values.

\subsection{Slice Sampling}

In this section we describe the developments in slice sampling and then show how to apply them to our optimisation problem. 
Suppose that we wish to sample from a possibly high dimensional un-normalised density
$ \pi(\bx) $. We do this by sampling uniformly from the region that lies under the density plot of $\pi$. This idea is formalised
by letting $u$ be an auxiliary ``slice-variable'' and defining a joint distribution $ \pi(\bx,u)$ that is uniform
on the set $ U = \{ ( \bx , u) : 0 < u < \pi(\bx) \} $.  Therefore, $ p ( \bx , u) = 1 / Z $ on $U$ and zero otherwise.
Here $ Z = \int_{\mathcal{X}} \pi(\bx ) d \bx $ is the appropriate normalisation constant.
The marginal is the desired normalised density as
$$
\pi(\bx) = \int_U \pi(\bx,u) d u = (1 / Z ) \int_0^{\pi(\bx)} d U = \pi(\bx)/Z  \; .
$$
We are then left with sampling from the uniform density on $U$. Neal (2003) provides a general slice algorithm. When it is straightforward
to sample from the ``slice'' region defined by $u$, namely
$ \mathcal{S}_u = \{ \bx : u < \pi(\bx) \} $, then a simple Gibbs sampler which iterates between drawing a uniform $ (u|\bx) \sim Uni ( 0, \pi(\bx) )  $
and $ (\bx|u) \sim Uni ( \mathcal{S}_u  ) $ provides a Markov chain with the joint distribution $\pi(\bx,u)$ and hence we can obtain marginal draws from
$ \pi(\bx)/Z$.

The Swenden-Wang algorithm (Edwards and Sokal, 1988) extends the application of slice sampling to product functionals.
Suppose that we wish to sample from a density that is a product of functions: 
$$
\pi_K(\bx) = \pi_1 (x ) \ldots \pi_K (x) / Z_K = \prod_{i=1}^K \pi_i(\bx) / Z_K \; .
$$ 
To do this, we introduce $K$ auxiliary uniform slice variables $ (u_1 , \ldots , u_K )$ with
a joint $\pi (\bx, u_1 , \ldots, u_K )$ defined to be uniform on the ``slice'' region:
$$
 \mathcal{S} = \{ (\bx,u) : u_i <  \pi_i(\bx) \; , \; 1 \leq i \leq K \; \}.
$$
The marginal distribution $\pi_K(\bx) = \int \pi(\bx ,u) d u = \prod_{i=1}^K \pi_i(\bx) / Z_K $.
We can sample this distribution using the complete conditionals
$$
(u_i | \bx ) \sim Uni ( 0, \pi_i(\bx) ) \; {\rm for} \;  i = 1 , \ldots , K  \; \; {\rm  and} \;  \; (\bx|u) \sim Uni ( \mathcal{S}_u ) \; .
$$
where $ \mathcal{S}_u = \{ \bx : u_i <  \pi_i(\bx) \; , \; 1 \leq i \leq K \; \} $.

We will be interested in additive objective functions:
$ f(\bx) = \sum_{i=1}^K f_i (\bx) $. Define $ \pi_i (\bx) = \exp ( - \kappa f_i (\bx) ) $. Now the Boltzmann distribution is
$$
\pi_K ( \bx ) = \exp \left ( - \kappa f (\bx) \right )/ Z_\kappa  = \exp \left ( - \kappa \sum_{i=1}^K f_i (\bx) \right ) / Z_K \;.
$$
The exponential slice sampler extends slice sampling to additive functionals by
letting $ y = ( y_1 , \ldots , y_K ) $ be a vector of exponential variables with
$ (y_i |\kappa ) \sim Exp(\kappa)$. The joint distribution is given by:
$$
\pi(  \bx , y_1 , \ldots , y_K ) = \exp \left ( - \kappa \sum_{i=1}^K y_i \right ) \prod_{i=1}^K \mathbb{I}\left ( 0 \leq y_i \leq f_i(\bx) \right ) / Z_\kappa.
$$
The MCMC algorithm uses 
the complete conditional  $ p( \bx | y_1 , \ldots , y_K )$ and $ p( y_i | \bx ) $ where 
$$
\pi( \bx|y )  \sim Uni \left \{   \{ f_i(\bx) \geq y_i \; , \; \forall i \; \}  = \cup_{i=1}^K \left ( \bx_i \in f^{-1}_i ( y_i ) \right ) \right . \}.
$$
The auxiliary variables are sampled from $ ( y_i | \bx ) \sim Exp( \kappa ) \mathbb{I} ( 0 , f_i (\bx) )$.

So far we have constructed a Markov chain on the augmented space $(\bx,y)$ which converges in distribution to $\pi(\bx,y)$:
we write $ \left ( \bx^{(n)} , y_1^{(n)} , \ldots , y_K^{(n)} \right ) \stackrel{D}{=} ( \bx , y_1,\ldots ,y_K ) \sim \pi(\bx,y) $ 
as $ n \rightarrow \infty$.
Given weak convergence, we also have for any functional that
$$
F \left ( \bx^{(n)} , y_1^{(n)} , \ldots , y_K^{(n)} \right ) \stackrel{D}{=} F ( \bx , y_1,\ldots , y_K ) 
$$
Hence we can project that joint draws back to the original space and view $ \bx^{(n)} $ as traversing the Boltzmann funcion
or equivalently $f(\bx)$.
This allows us to traverse the modes of interest and has the
advantage of scalability in $K$.

A related approach is to allow the energy level, $\kappa$, to be random. One 
places pseudo-prior weights, denoted by $p(\kappa)$, and simulates the mixture Boltzmann distribution $ \sum_{\kappa} p(\kappa) \pi_\kappa (\bx) $.
Another line of research is based on simulated annealing (Kirkpatrick et al, 1983, Aarts and Korst, 1988, Van Laarhoven and Aarts, 1987)
which increases the energy level with the length of simulation in an appropriate schedule (Gidas, 1985).
Other approaches that randomize $\kappa$ include multi-canonical sampling (Berg and Neuhaus, 1982),
simulated tempering (Marinari and Parisi, 1992, Geyer and Thompson, 1995) uses a random walk on a set of
energy level, equi-energy sampling (Kou, Zhou and Wong, 2006), evolutionary
MCMC (Liu, Liang and Wong, 2000, 2001) and the Wang-Landau algorithm (Liang, 2005, Atchade and Liu, 2010) and 
When $ \mathcal{X}$ is bounded and $ supp (f) < \infty $, slice sampler has the added property of
geometric convergence (Roberts and Polson, 1994) and uniformity (Mira and Tierney, 2002).

\section{Four Examples}

Figures 1 and 2 show contour and drape-plots of the Rosenbrock, Himmelblau, Rastigrin and Shubert functions.
This set of functions exhibits a variety of challenges: a global mode in a long valley (Rosenbrock),
multiple local and one global mode (Rastrigin) to multi-modality (Himmelblau, Shubert). For example, 
traditional derivative-based methods have difficulties for the Rosenbrock function to traverse its long steep valley.

Let $f(\bx)= \sum_{i=1}^K f_i(x_1,x_2)$ be a  bivariate additive objective function defined over a bounded region for 
$ \bx = \{ x_1 , x_2 \} \in \Re^2 $.
We will apply the exponential slice sampler as the  functions to the Boltzmann distribution
$$
\pi_K ( x_1 , x_2 ) = \exp \left ( - \kappa \sum_{i=1}^K f_i(x_1,x_2) \right ) / Z_K \; .
$$
First, as in simulated tempering, we have to define a set of temperatures $ \{ \kappa_1 < \ldots < \kappa_m \} $ to run our Markov chain.
For all four examples, we pick $m=4$, and we use $ \kappa \in \{ 0.1, 0.5 , 1 , 5 \} $ except for the Rosenbrock function where
we set  $ \kappa \in \{ 1,5 ,50,5000 \} $ which requires higher energy levels. 

\subsection{Rosenbrock function}

Rosenbrock's valley is a classic optimization problem that illustrates the difficulties with local search methods.
The global minimum lays inside a long, narrow flat valley. Finding the valley is straightforward; however, getting to the minimum is hard.
We need to find the minimum $(x_1,x_2)=(1,1)$ of the function:
$$
f(x_1,x_2) = ( 1- x_1 )^2 + c(x_2-x^2_1)^2 \; \; {\rm where} \; c=100 \; .
$$
The Boltzmann distribution has density
$$
\pi_K ( x_1 ,x_2) = \exp \left ( - \kappa \left \{ ( 1- x_1 )^2 + c(x_2-x^2_1)^2 \right \} \right ) / Z_\kappa \; .
$$
There are a number of ways of introducing slice variables. We choose to slice out the last nonlinear factor.  
Let $ ( u|x_1,x_2) \sim Uni \left (0, \exp \left \{ -  \kappa c (x_2-x^2_1)^2 \right \} \right ) $. 
Then we have a three variable joint distribution:
$$
\pi_{\kappa } ( x ,y , u ) = \exp \left \{ - \kappa ( 1- x_1 )^2 \right \}
 \mathbb{I} \left ( 0 \leq u \leq \exp \left \{ - \kappa c (x_2-x^2_1)^2 \right \} \right ) / Z_\kappa .
$$
We can implement MCMC using a partially collapsed Gibbs sampler (van Dyk and Park, 2008, Park and van Dyk, 2009). 
That is, we can marginalise $u$ out of the draw for $ x_2$
and use the conditional $ \pi( x_2 | x_1 ) $ rather than $ \pi( x_2 | x_1 , u )$.
The complete conditionals are then:
\begin{align*}
\pi ( x_2 | x_1 ) & \sim \mathcal{N} ( x^2_1 , 2 /  \kappa c ),\\
  \pi ( x_1 | x_2,u ) & \sim \mathcal{N} ( 1 , 2/ \kappa  ) \mathbb{I} \left ( a(u,x_2) \leq x \leq b(u,x_2) \right ),\\
\pi ( u | x_1 , x_2 ) & \sim Uni \left (0, \exp \left \{ - \kappa c (x_2-x^2_1)^2 \right \} \right ).
\end{align*}
The interval $(a(u,x_2), b(u,x_2)) $ is found by inverting the slice region:
$$
u \leq \exp \left \{ - \kappa c (x_2-x^2_1)^2 \right \} \; {\rm implies} \; x_2 - \sqrt{ - \ln u / \kappa c} \leq x^2_1 \leq  x_2 
+ \sqrt{ - \ln u / \kappa c}.
$$
For $b>0$ and $ a \leq x^2_1 \leq b$, we have $ - \sqrt{a} \leq x_1 \leq \sqrt{b} $ and so
$$
a( u, x_2 ) = - \sqrt{ x_2 - \sqrt{ - \ln u / \kappa c} } \; {\rm and} \; b( u, x_2 ) = \sqrt{ x_2 + \sqrt{ - \ln u / \kappa c} }.
$$
Figure 3 shows a sensitivity analysis for a range of energy values $ \kappa \in \{ 1,5 ,50,5000 \} $. We run our MCMC algorithm
for $G=1000$ with a burn-in of $G_0=100$. Higher energy levels are required for the chain to travel 
along the valley to the minimum at $(x_1,x_2)=(1,1)$. As we increase $\kappa $, the slice sampler is able to traverse the valley and find the minimum.

\subsection{Himmelblau's function}

Himmelblau's function is defined by 
$$
f(x_1,x_2) =  ( x^2_1 + x_2 - 11 )^2 + (  x_1 + x_2^2 - 5 )^2.
$$
This function has four identical local minima at zero and a local maximum at $ (x_1, x_2) = (-0.27,-0.92) $.
The minima are at 
$$ 
(3,2) , (-2.805, 3.131 ), (-3.779 , -3.282) , (3.584 , -1.848 ) 
$$ 
with a function value of zero. The associated Boltzmann distribution is 
$$
\pi_\kappa ( x_1 , x_2 ) = \exp \left \{ - \kappa \left (  ( x^2_1 + x_2 - 11 )^2 + (  x_1 + x_2^2 - 5 )^2 \right ) \right \} / Z_\kappa \; .
$$
Due to the quadratic terms also containing squares this distribution is not straightforward to simulate from.
We observe that the following inequalities hold for the minima: $ x_1^2 + x_2 - 11 <0 $ and $  x_1 + x_2^2 - 5 > 0 $.

To implement slice sampling, we use a latent variable augmentation $ u = ( u_1, u_2 ) $ and a joint distribution 
$$
\pi_\kappa ( x_1 , x_2 , u_1 , u_2 ) = \mathbb{I} \left ( 0 \leq u_1 \leq \exp \left \{ - \kappa ( x^2_1 + x_2 - 11 )^2 \right \} \right )
 \mathbb{I} \left ( 0 \leq u_2 \leq \exp \left \{ - \kappa (  x_1 + x_2^2 - 5 )^2  \right \} \right ) / Z_\kappa \; .
$$
Given $(u_1,u_2)$, we can invert the slice regions to obtain the inequalities:
$$
-  \kappa^{-1} \log u_1 \geq  ( x^2_1 + x_2 - 11 )^2 \; \; {\rm and} \; \; -  \kappa^{-1} \log u_2 \geq  ( x_1 + x_2^2 - 5 )^2.
$$
Therefore, for $(x_1|x_2)$, we have
\begin{align*}
a_1 = 11 - x_2 - \sqrt{ -  \kappa^{-1} \log u_1 } & \leq x^2_1 \leq  11 - x_2 + \sqrt{ -  \kappa^{-1} \log u_1 } =b_1;\\
c_1 = 5 - x_2^2 - \sqrt{ -  \kappa^{-1} \log u_2 }& \leq x_1 \leq  5 - x_2^2 + \sqrt{ -  \kappa^{-1} \log u_2 } =d_1.
\end{align*}
Given the inequalities: $ a_1 \le x_1^2 \leq b_1 $ and $ c_1 \leq x_1 \leq d_1 $,  we can first without loss of generality  replace $a_1 \rightarrow \max(a_1,0) $ and
assume that $ a_1 \geq 0 $.  Then, we have the union of the following regions:
$$
- \sqrt{b_1} \leq x_1 \leq - \sqrt{a_1} \; {\rm and} \; \sqrt{a_1} \leq x_1 \leq b_1 .
$$
Combining, we have 
$$ 
\max \left (  - \sqrt{ b_1 } ,c_1 \right ) \leq x_1 \leq \min \left ( - \sqrt{a_1} , d_1 \right ) \; {\rm or} \;
\max \left (   \sqrt{ a_1 } ,c_1 \right ) \leq x_1 \leq \min \left (  \sqrt{b_1} , d_1 \right )   
$$
For the conditional $\pi (x_2|x_1)$ we can argue in a similar fashion to obtain the constraints
\begin{align*}
a_2 = 5 - x_1 - \sqrt{ -  \kappa^{-1} \log u_1 } & \leq x_2^2 \leq  5 - x_1 + \sqrt{ -  \kappa^{-1} \log u_1 } =b_2;\\
c_2 = 11 - x_1^2 - \sqrt{ -  \kappa^{-1} \log u_2 }& \leq x_2 \leq  11 - x_1^2 + \sqrt{ -  \kappa^{-1} \log u_2 } =d_2.
\end{align*}
The complete set of conditionals is then given by:
\begin{align*}
\pi( x_1| x_2 , u_1 , u_2 ) & \sim Uni \left ( \max \left (  - \sqrt{ | a_1| } ,c_1 \right ) , \min \left ( \sqrt{b_1} , d_1\right ) \right ),\\
\pi( x_2| x_1 , u_1 , u_2 ) & \sim Uni \left ( \max \left (  - \sqrt{ | a_2| } ,c_2 \right ) , \min \left ( \sqrt{b_2} , d_2 \right ) \right ),\\
\pi( u_1 | x_1, x_2 ) & \sim Uni \left ( 0 , \exp \left \{ - \kappa ( x^2_1 + x_2 - 11 )^2 \right \} \right ), \\
\pi( u_2 | x_1 , x_2) & \sim Uni \left ( 0 , \exp \left \{ - \kappa ( x_1 + x_2^2 - 5 )^2 \right \} \right ).
\end{align*}
Figure 4 shows a sensitivity analysis to $ \kappa \in \{ 0.1,0.5 ,1,5 \} $. The slice sampler is again run for only $G=1000$ iterations
with a burn-in of $G_0=100$. With longer chains and larger energy levels the algorithm will traverse the four modes with equal probability.
Of the examples that we consider here, the Himmleblau function would benefit the most from a mixture energy level distribution to traverse the contours of the
underlying function.

\subsection{Rastrigin}

The $2$-dimensional Rastrigin function is defined on the region $ -5.12 < x_j < 5.12 $ by:
$$
f( x_1, x_2 ) = 2 A + \sum_{j=1}^2 \left ( x_j^2 - A \cos ( 2 \pi x_j ) \right ) \; \; {\rm with} \; \; A=10 
$$
with a global minimum at $ (x_1, x_2)=(0,0)$.
The Boltzmann distribution then becomes 
$$
\pi_\kappa ( x_1, x_2  ) = \exp \left \{ - \kappa \sum_{j=1}^2 x_j^2 \right \}
 \exp \left \{ \kappa A \sum_{j=1}^2  \cos ( 2 \pi x_j ) \right \} / Z_\kappa \; .
$$
We use exponential slice variables $(y_1,y_2)$ and a joint distribution defined by  
$$
 \pi_\kappa ( x_1 , x_2 , y_1 , y_2  ) = \exp \left \{ - \kappa \sum_{j=1}^2 x_j^2 \right \}
\mathbb{I} \left ( - \kappa A \cos (2 \pi x_j ) \leq y_j  \right ) e^{- y_j} / Z_\kappa \; .
$$
The slice region is invertible and is specified by the set of inequalities for $j=1,2$ 
$$
 \cos (2 \pi x_j) \geq  \left ( - y_j /A \kappa \right ) \; .
$$
The conditional $x_j$ draw then results from a normal draw restricted to this interval set.
A Gibbs sampler with the conditionals for $ j=1,2 $ is as follows:
\begin{align*}
\pi_\kappa ( x_j | x_{-j} , y ) & \sim \mathcal{N} \left ( 0 , (2\kappa)^{-1} \right )
 \mathbb{I} \left ( \cos (2 \pi x_j) \geq  \left ( - y_j /A \kappa \right )  \right )\\
\pi_\kappa ( y_j | x_j  ) & \sim Exp [ - \kappa A \cos (2 \pi x_j ) , \infty )
\end{align*}
For the implementation over $x_j\in [-5.12,5.12]$, we draw the truncated normals and truncated exponentials for the slice variables.

Figure 5 shows a sensitivity analysis to $ \kappa \in \{ 0.1,0.5 ,1,5 \} $ with $G=1000$ and a burn-in of $G_0=100$.
Again slice sampling of the Boltzmann distribution finds the mode in a straightforward manner.

\subsection{Shubert function}

The Shubert function is defined  within the domain $ \mathbb{I} (-10, 10) $ by
$$
f(x_1,x_2)=-C(x_1)C(x_2)
\; \; {\rm where} \; \;  C( x ) = \sum_{j=1}^5 j \cos \left ( ( j+1 ) x + j \right ) \; .
$$ 
We need to simulate from the Boltzmann distribution
$$
\pi_\kappa ( x_1 , x_2 ) = e^{ - \kappa  \sum_{j=1}^5 j \cos \left ( ( j+1 ) x_1 + j \right ) \cdot \sum_{j=1}^5 j \cos \left ( ( j+1 ) x_2 + j \right )} / Z_\kappa
$$
For the conditional $ \pi_\kappa ( x_1 | x_{2})$, we can write
$$
 \pi_\kappa ( x_1 | x_2 ) = \prod_{j=1}^5 e^{ \kappa C(x_2) j \cos \left ( (j+1) x_1 + j \right ) }.
$$
We introduce a set of auxiliary slice variables
$ y_{j} , 1 \leq j \leq 5 $ for each $x_i, i=1,2$ that are  conditional exponentials. The corresponding joint is:
$$
 \pi_\kappa ( x_1 , x_2 , y_{1} , \ldots , y_{5} ) = \prod_{j=1}^5 \exp (-\kappa y_j)\mathbb{I}
 \left ( y_{j} \geq C( x_{2} ) j \cos \left ( (j+1) x_1 + j \right ) ) \right ) / Z_\kappa
$$
This is inverted using the condition:
$$
 \left (   \frac{y_j}{jC(x_{2})} \right )\mathbb{I}_{C(x_2)<0}\le  \cos( (j+1) x_i + j )\leq  \left (   \frac{y_j}{jC(x_{2})} \right )\mathbb{I}_{C(x_2)>0} \; .
$$
This gives a collection of intervals $ I_1(y_j,x_2)$ for $ 1 \leq j \leq 5 $
for each of the slice variables; $x_1$ is then uniformly distributed over the intersection of these intervals, $\cap_j I_1(u_j,x_2)$.

Hence we can then run a Gibbs sampler with the conditionals, for $ 1 \leq j \leq 5 $:
\begin{align*}
\pi_\kappa ( x_1 |  x_2 , y ) & \sim Uni (\cap_{j=1}^5 I_1(y_j,x_2)), \\
\pi_\kappa ( y_j | x_1  ) & \sim Exp \left [  C( x_2 ) j \cos \left ( (j+1) x_1 + j \right ), \infty  \right ) .
\end{align*}
 Similarly, this process is repeated for $ \pi_\kappa ( x_2 | x_1 ) $. This defines a $12$-dimensional Gibbs sampler that is able
to traverse the joint distribution.

Figure 6 shows a sensitivity analysis to the same set of energy levels used for the Rastgrin and Himmleblau functions,
namely $ \kappa \in \{ 0.1,0.5 ,1,5 \} $. Again the projected draws from the uniform slice region traverse the modes of the associated
Boltzmann distribution in an efficient manner.

\section{Discussion}

We have described how slice sampling methods can be applied to functional optimisation. Our approach is parallelisable as in slice sampling (Tibbits et al, 2011).
While we have only considered four test
functions our methodology applies to a multitude of multi-modal functions, see Molga and Smutnicki (2005) for a list of possible candidates.
Our approach is flexible enough to handle additive functions that are multiplicative. For example, the 
Michalewicz function is defined by 
$f(x_1,x_2) = - \sum_{i=1}^2 sin(x_i) sin^{2m} \left ( i x_i^2 / \pi \right ) $
with $m=10$ (see Yang, 2010a,b). This function has a minimum of $ - 1.801 $ at the point $(x_1,x_2)=( 2.20319,1.57049)$ and it is challenge to fix the minimum.
We also note that certain functions are straightforward as the Boltzmann distribution is conditionally normal.
For example, the Booth function defined by
$f(x_1,x_2) = ( x_1 + 2 x_2 - 7 )^2 + ( 2 x_1 + x_2 - 5 )^2 $ has a minimum of zero at the point $(x_1,x_2)=(1,3)$.
The minimum can be identified without resorting to simulation as the Booth function can be factorised into the quadratic form 
$ ( \bx - \mu )^\prime Q^{-1} ( \bx - \mu ) $
where $ x = (x_1 , x_2 ) $ with $ \mu = ( 1 , 3 ) $ and $ Q = \frac{1}{9} \left ( 5 , 4 ; - 4 , 5 \right ) $.
Therefore, the minimum can be immediately identified without resorting to simulation.

There are clearly many other applications of these methods. For example, simulated annealing methods have been proposed in mixed integer non-linear programming problems
(Cardoso et al, 1997) and in constrained optimisation (Geman and Geman, 1985, 
Geman, 1990, Whittle, 1992, Birge and Louveaux, 1997, Mueller, 2000, Asmussen and Glynn, 2008). Slice sampling the Boltzmann distribution provides a
flexible alternative. 

\section{References}

\begin{description}

\item Aarts, E and Korst, J. (1989). \textit{Simulated Annealing and Boltzmann Machines}. Wiley, NY.

\item Atchade, Y. and J. Liu (2010). The Wang-Landau algorithm in general state spaces: applications and convergence analysis.
\textit{Statistica Sinica}, 20, 209-233.

\item Asmussen, S. and P. Glynn (2008). \textit{Stochastic Simulation}. Springer-Verlag, New York.

\item Azencott, R. (1988). Simulated Annealing. \textit{Seminaire Bourbaki}, 697.

\item Berg, B.A. and T. Neuhaus (1992). Multicanonical ensemble: A new approach to simulate first-order phase transitions.
\textit{Physical Review Letters}, 68, 9-12.

\item Birge, J.R. and F. Louveaux (1997). \textit{Introduction to Stochastic
Programming}. Springer, New York.

\item Cardoso, M.F., R.L. Salcedo, S. Feyo de Azevedo and D. Barbosa (1997). A simulated annealing approach to the solution of mnlp problems.
\textit{Computers Chem. Engng}, 21(12), 1349-1384.\medskip

\item Devroye, L. (1986). \textit{Non-Uniform Random Variate Generation}, Springer Verlag, New York.

\item Edwards, R.G. and A.D. Sokal (1988). Generalisation of the Fortuin-Kastelyn-Swendsen-Wang algorithm.
\textit{Phys. Review D}, 38(6), 2009-2012.

\item Geman, D. (1990).
\textit{Random Fields and Inverse Problems in Imaging}. Lecture Notes, Springer-Verlag, 113-193.

\item Geman, D. and S. Geman (1985). Relaxation and annealing with constraints.
\textit{Technical Report 35}, Brown University.

\item Geyer, C.J. (1992). Practical Markov chain Monte Carlo
(with discussion). \textit{Statistical Science}, 7, 473-511.

\item Geyer, C.J. and E.A. Thompson (1995). Annealing MCMC with applications to Ancestral Inference.
\textit{Journal of the American Statistical Association}, 90, 909-920.

\item Gidas, B. (1985). Nonstationary Markov chains and convergence of the annealing algorithm.
\textit{J. Stat. Phys.}, 39, 73-131.

\item Gramacy, R.B. and H. Lee  (2011). Optimisation under unknown constraints.
\textit{Bayesian Statistics}, 9, 229-257.

\item Gramacy, R. and M. Taddy (2010). Categorical Inputs, Sensitivity Analysis, Optimization and Importance Tempering in {\tt tgp}.
\textit{J. Statistical Software}, 22(6), 1-48.

\item Janson, S. and M. Middendorf (2005). A Hierarchical Particle Swarm Optimiser and its Adaptive Variant.
\textit{IEEE Trans in Systems, Man, and Cybernetics}, 35(6), 1272-1282.

\item Kennedy, J. and R. Eberhart (1995). Particle Swarm Optimisation.
\textit{IEEE Int. Conf. on Neural Networks}, 1942-1948.

\item Kirkpatrick, S., C.D. Gelatt and M.P. Vecchi, (1983). Optimization by
simulated annealing, \textit{Science}, 220, 671-680.

\item Kou, S.C., Zhou, Q. and W.H. Wong (2006). Equi-Energy Sampler with applicationbs in statistical inference
and statistical mechanics (with Discussion). \textit{Annals of Statistics}, 34(4), 1581-1619.

\item Leman, S.C., Y. Chen and M. Lavine (2009). The Multi-Set Sampler.
\textit{Journal of the American Statistical Association}, 104, 1029-1041.

\item Liang, F. (2005). Generalized Wang-Landau algorithm for Monte Carlo computation.
\textit{Journal of American Statistical Association}, 100, 1311-1337.\medskip

\item Liang, F., C. Liu and R.J. Carroll (2007). Stochastic approximation in Monte Carlo computation.
\textit{Journal of American Statistical Association}, 102, 305-320.

\item Liu, J.S., Liang, F. and W.H. Wong (2000). The Multiple-Try Method and Local Optimisation in
Metropolis Sampling. \textit{Journal of the American Statistical Association}, 95, 121-134.

\item Liu, J.S., Liang, F. and W.H. Wong (2001). A theory of dynamic weighting in Monte Carlo.
\textit{Journal of the American Statistical Association}, 96, 561-573.

\item Marinari, E. and G. Parisi (1992). Simulated Tempering: A Monte Carlo scheme.
\textit{Euro Phys. Lett. EPL}, 19, 451-458.

\item Mira, A. and L. Tierney (2002). Efficiency and Convergence Properties of Slice Samplers.
\textit{Scandinavian Journal of Statistics}, 29(1), 1-12.

\item Molga, M. and C. Smutnicki (2005). Test functions for optimization needs. \textit{Working Paper}.

\item Mueller, P. (2000). Simulation-Based Optimal Design.
\textit{Bayesian Statistics}, 6, 459-474.

\item Neal, R. (2003). Slice Sampling (with Discussion). \textit{Annals of Statistics}, 31(3), 705-767.

\item Park, T. and D.A. van Dyk (2008). Partially Collapsed Gibbs Samplers: Theory and Methods.
\textit{Journal of the American Statistical Association}, 103, 790-796.

\item Pincus, M. (1968). A Closed Form Solution of Certain Dynamic
Programming Problems. \textit{Operations Research}, 16, 690-694.

\item Pincus, M. (1970). A Monte Carlo Method for Approximate Solution
of certain types of Constrained Optimization Problems.
\textit{Operations Research}, 18(6), 1225-1228.

\item Polson, N.G. (1992). Comment on ``Practical Markov chain
Monte Carlo''. \textit{Statistical Science}, 7, 490-491.

\item Polson, N. G. (1996). Convergence of Markov Chain Monte Carlo Algorithms.
\textit{Bayesian Statistics 5}, Bernardo et al eds, Oxford, 297-321.

\item Polson, N. G., J. G. Scott and J. Windle (2012). The Bayesian Bridge. \textit{Working Paper}.

\item Roberts, G.O.  and Polson, N. G. (1994). On the Geometric Convergence of the Gibbs Sampler.
\textit{J. R. Statist. Soc.}, B, 56(2), 377-384.

\item Roberts, G.O.  and J. Rosenthal (1999). Convergence of slice sampler Markov chains.
\textit{J. R. Statist. Soc.}, B, 61, 643-660.

\item Taddy, M., H.K.H. Lee, G.A. Gray and J.D. Griffin (2009). Bayesian Guided Pattern Search for Robust Local Optimization.
\textit{Technometrics}, 51, 389-401.

\item Tibbits, M., M. Haran, J.C. Lietchy (2011). Parallel Multivariate Slice Sampling.
\textit{Technical Report}.

\item Tierney, L. (1994). Markov Chains for exploring Posterior
Distributions (with discussion). \textit{Annals of Statistics}, 22, 1701-1786.

\item Tweedie, R and K. Mengersen (1994). Rates of convergence of the Hastings-Metropolis algorithms.
\textit{Annals of Statistics}, 24, 101-121.

\item van Dyk, D.A. and T. Park (2009). Partially Collapsed Gibbs Samplers: Illustrations and Applications.
\textit{Journal of Computational and Graphical Statistics}, 18, 283-305.

\item Van Laarhoven P.J. and Aarts, E.H.J. (1987). \textit{Simulated Annealing:
Theory and Applications}, CWI Tract 51, Reidel, Amsterdam.

\item Wang, F. and D.P. Landau (2001). Efficient multi-range random walk algorithm to calculate the density of states.
\textit{Phys. Rev. Lett}, 86, 2050-2053. 

\item Whittle, P. (1982). \textit{Optimization over Time} (Volume 1). Wiley.\medskip

\item Yang, X-S (2010a). Test problems in Optimization. \textit{Working Paper}, University of Cambridge.\medskip

\item Yang, X-S (2010b). Firefly algorithms for Multimodal Optimization. \textit{Working Paper}, University of Cambridge.\medskip

\item Zhou, Q. and W.H. Wong (2008). Reconstructing the energy landscape of a distribution of Monte Carlo samples.
\textit{Annals of Applied Statistics}, 2(4),1307-1331. \medskip

\end{description}

\newpage
\begin{figure}[hbp]
\vspace{-0.2in}
\includegraphics[height=6in,width=\textwidth]{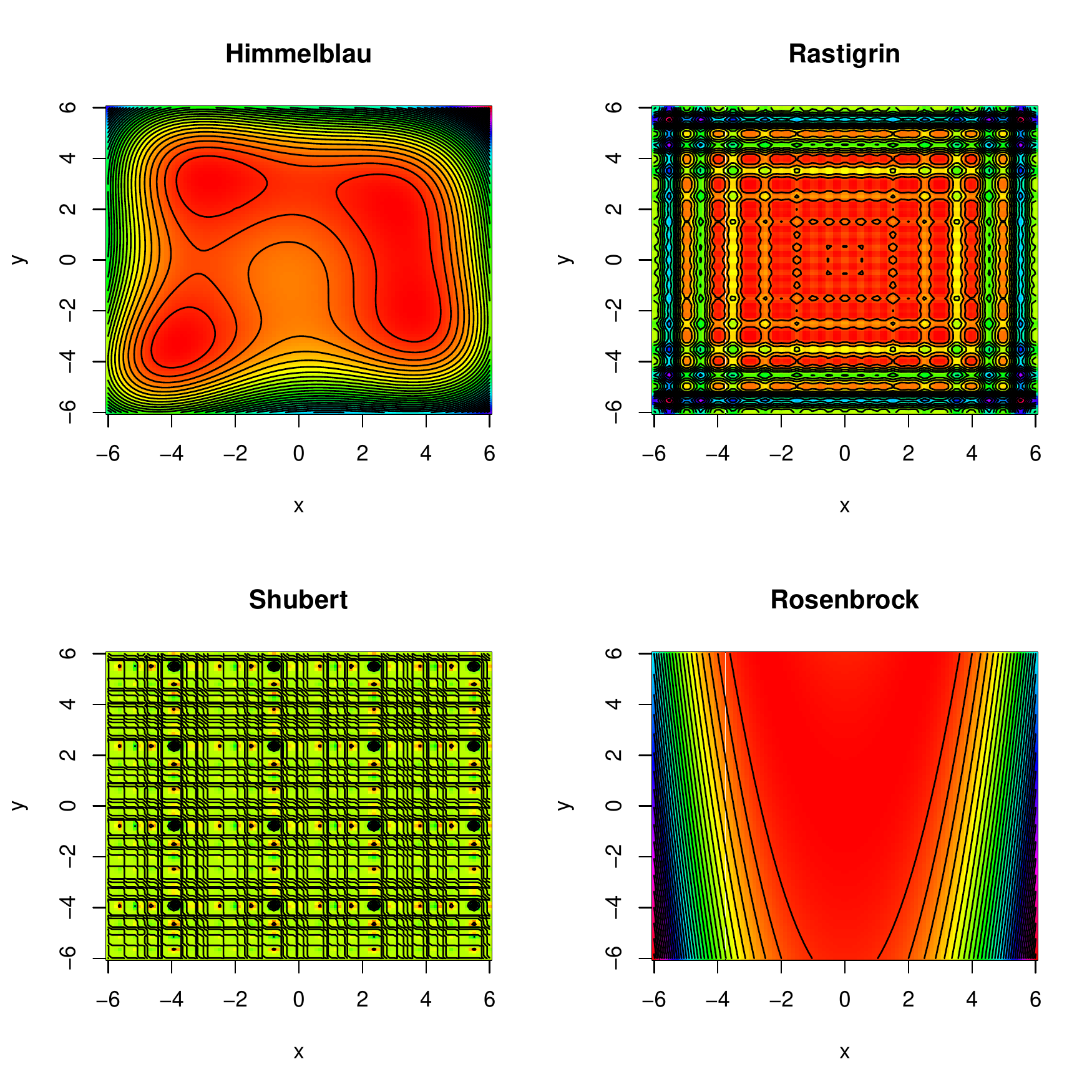}
\caption{Contour plots: Himmelblau, Rastigrin, Shubert and Rosenbrock}
\end{figure}
\newpage
\begin{figure}[hbp]
\vspace{-0.2in}
\includegraphics[height=6in,width=\textwidth]{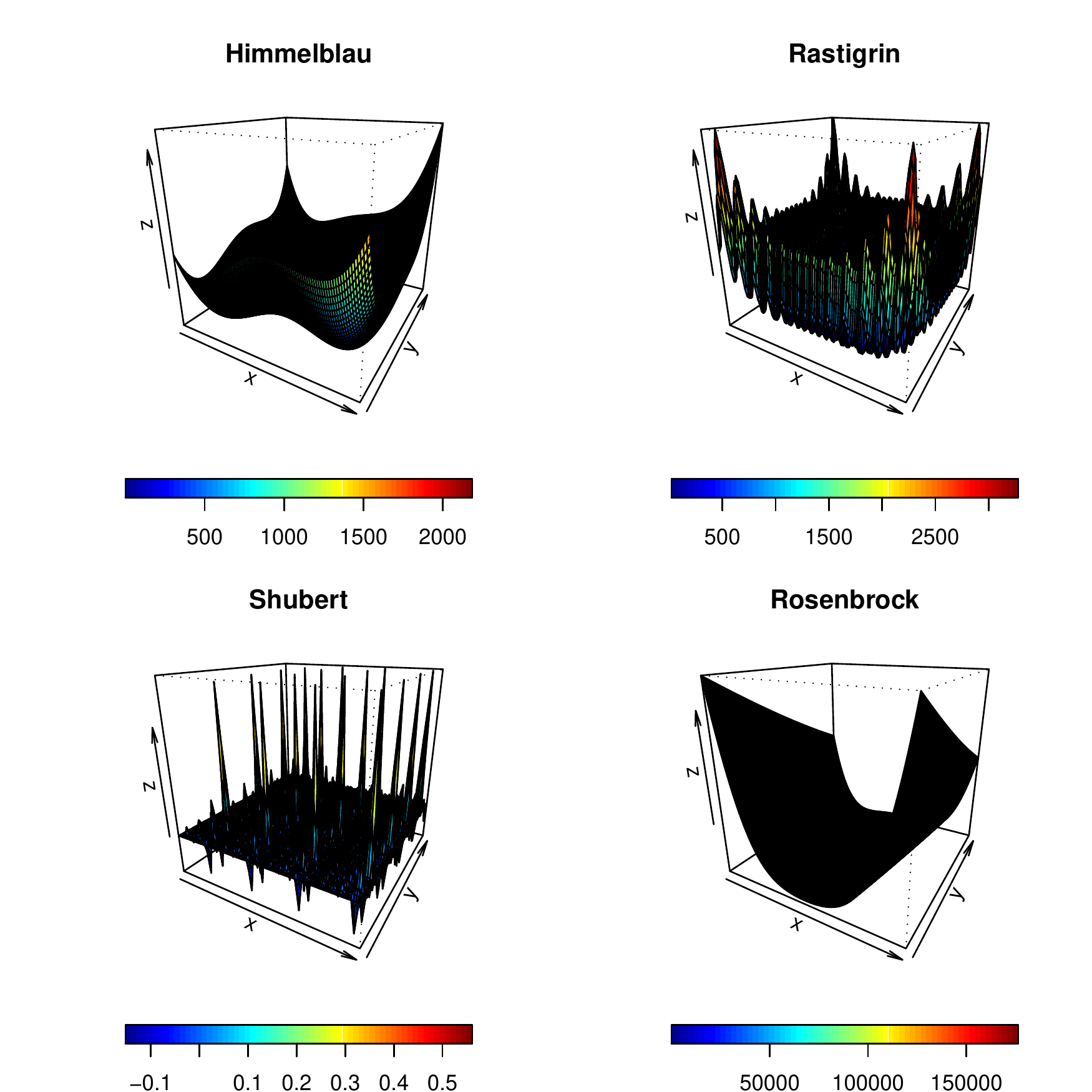}
\caption{Drape plots: Himmelblau, Rastigrin, Shubert and Rosenbrock}
\end{figure}
\newpage
\begin{figure}[hbp]
\vspace{-0.2in}
\includegraphics[height=6in,width=\textwidth]{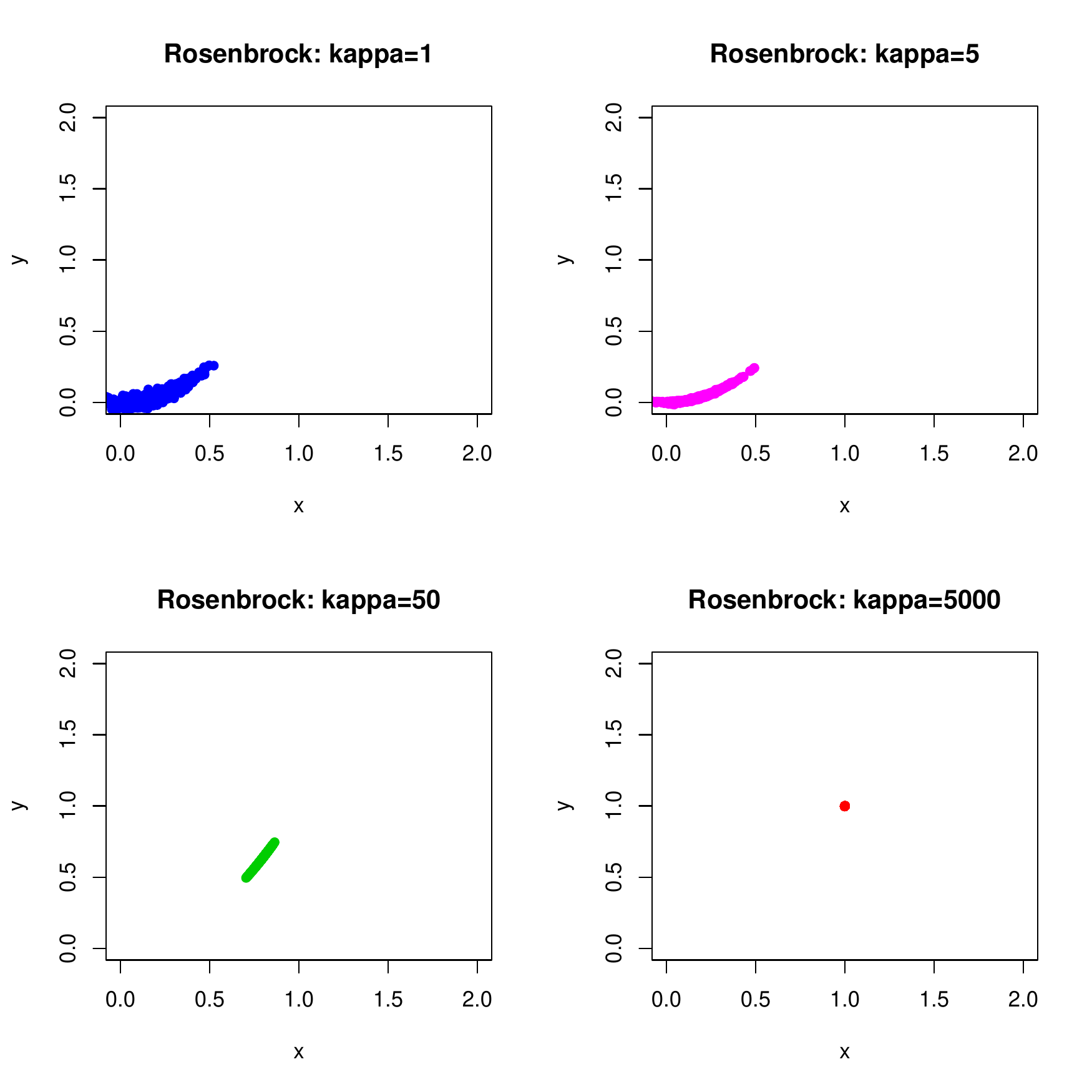}
\caption{Rosenbrock function: $f(x_1,x_2) = (1-x_1)^2 + 100 (x_2 -x_1^2)^2 $. The minimum occurs at $(x_1,x_2)=(1,1)$.}
\end{figure}
\newpage
\begin{figure}[hbp]
\vspace{-0.2in}
\includegraphics[height=6in,width=\textwidth]{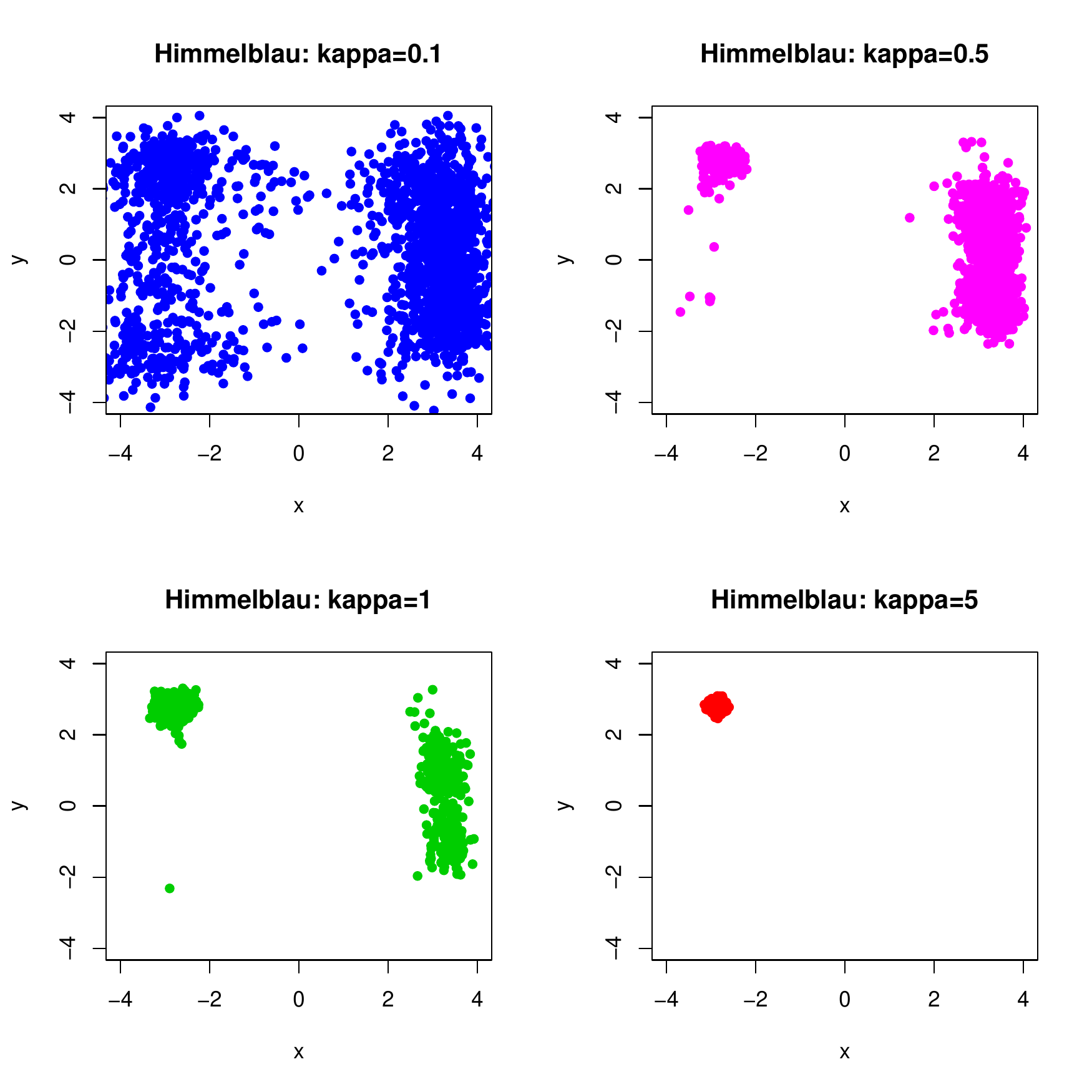}
\caption{Himmelblau function: $f(x_1,x_2) = ( x_1^2 +x_2 - 11 )^2 + ( x_1 + x_2^2 - 5)^2 $. There are four idenital local minima at zero and a local maximum at 
$ x_1 = -0.27 , x_2 = -0.92 $.}
\end{figure}
\newpage
\begin{figure}[hbp]
\vspace{-0.2in}
\includegraphics[height=6in,width=\textwidth]{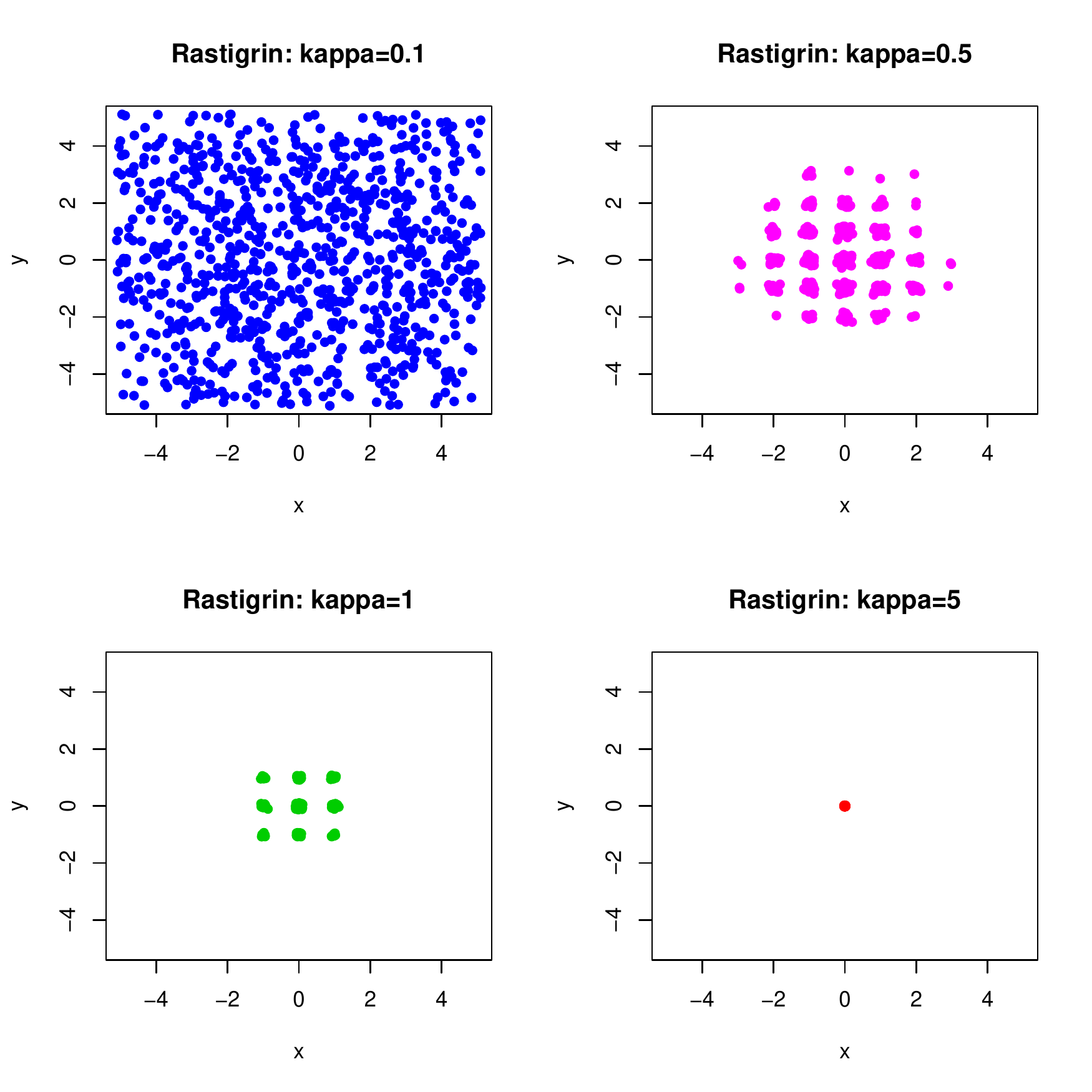}
\caption{$2$-d Rastigrin function: $f(x_1,x_2) = 2 A + \sum_{j=1}^2 ( x_j^2 - A cos(2 \pi x_j ) $ with $A=10$ and $-5.12<x_j<5.12$.}
\end{figure}
\newpage
\begin{figure}[hbp]
\vspace{-0.2in}
\includegraphics[height=6in,width=\textwidth]{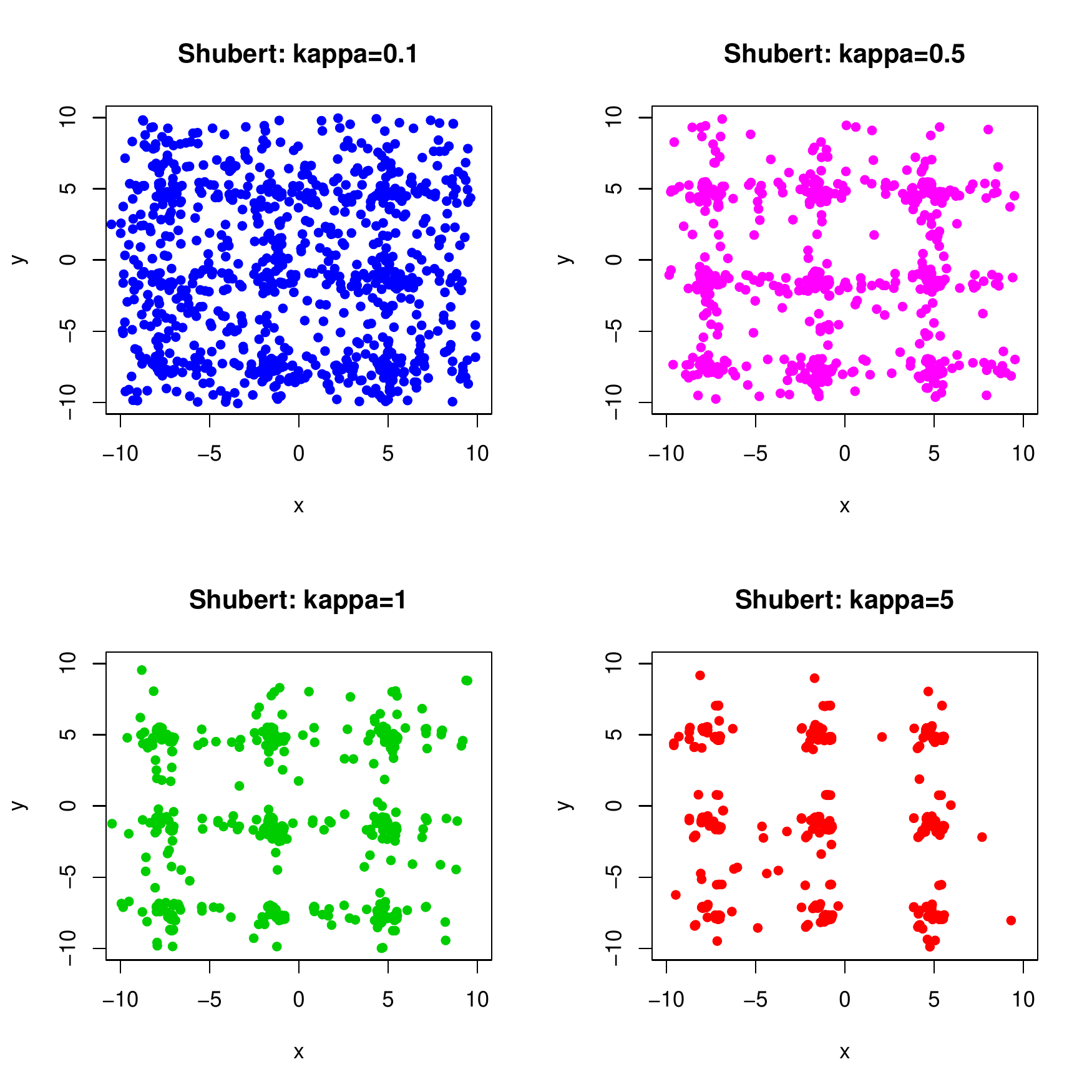}
\caption{Shubert function: $f(x_1,x_2) = - C(x_1 ) C(x_2) $ where $ C(x) = \sum_{j=1}^5 j cos \left ( (j+1) x + j \right ) $.}
\end{figure}

\end{document}